\documentclass{article}

\usepackage{amsmath,amsfonts,amssymb,a4,enumerate,epsfig,theorem,psfrag}
\usepackage[latin1]{inputenc}

\newcommand{\be}{{\phi_\pi}}

\newcommand{\sign}{{\operatorname{sign}}}
\newcommand{\signseq}{{\operatorname{signseq}}}
\newcommand{\trace}{{\operatorname{tr}}}

\newcommand{\transpose}{{\ast}}

\newcommand{\RR}{{\mathbb{R}}}

\newcommand{\Dd}{{\cal{D}}}
\newcommand{\Ee}{{\cal{E}}}

\newcommand{\Jj}{{\cal{J}}}

\newcommand{\Pp}{{\cal{P}}}

\newcommand{\Tt}{{\cal{T}}}
\newcommand{\Uu}{{\cal{U}}}

\newcommand{\bF}{{\mathbf{F}}}
\newcommand{\bL}{{\mathbf{L}}}
\newcommand{\bQ}{{\mathbf{Q}}}
\newcommand{\bR}{{\mathbf{R}}}
\newcommand{\bU}{{\mathbf{U}}}
\newcommand{\bT}{\tau}

\newcommand{\ILa}{\Tt_\Lambda}
\newcommand{\JLa}{\bar\Jj_\Lambda}
\newcommand{\JLao}{\Jj_\Lambda}
\newcommand{\PLa}{\Pp_\Lambda}
\newcommand{\UpLa}{\Uu^\pi_\Lambda}

\newcommand{\proof}{{\noindent\bf Proof: }}

\def\qed{\unskip\nobreak\hfil\penalty50\hskip1.75em\null\nobreak\hfil
$\blacksquare$ {\parfillskip=0pt \finalhyphendemerits=0 \par}\medbreak}

\newcommand\capsize{\relax}
\newcommand\nobf{\noindent\bf}

\newcommand\diag{\operatorname{diag}}

\newtheorem{lemma}{Lemma}[section]
\newtheorem{theo}[lemma]{Theorem}
\newtheorem{prop}[lemma]{Proposition}

\newtheorem{defin}[lemma]{Definition}

\title{An atlas for tridiagonal isospectral manifolds}
\author{Ricardo S. Leite, Nicolau C. Saldanha and Carlos Tomei }

\begin{document}
\maketitle

\begin{abstract}
Let $\ILa$ be the compact manifold of real symmetric tridiagonal matrices
conjugate to a given diagonal matrix $\Lambda$ with simple spectrum.
We introduce {\it bidiagonal coordinates},
charts defined on open dense domains forming an explicit atlas for $\ILa$.
In contrast to the standard inverse variables,
consisting of eigenvalues and norming constants,
every matrix in $\ILa$ now lies in the interior of some chart domain.
We provide examples of the convenience of these new coordinates
for the study of asymptotics of isospectral dynamics,
both for continuous and discrete time.
\end{abstract}

\medbreak

{\noindent\bf Keywords:} Jacobi matrices, tridiagonal matrices,
norming constants,
Toda flows,
$QR$ algorithm.

\smallbreak

{\noindent\bf MSC-class:} 65F18; 15A29.

\section{Introduction}

Let $\Lambda$ be a real diagonal matrix with simple spectrum and
$\ILa$ be the manifold of real, symmetric, tridiagonal matrices having
the same spectrum as $\Lambda$.
The purpose of this paper is to present an explicit atlas for $\ILa$:
the charts in the atlas define the \textit{bidiagonal coordinates}
on open dense subsets of $\ILa$.
As is familiar to numerical analysts, many algorithms
to compute spectra operate by iteration on Jacobi matrices,
yielding approximations of reduced tridiagonal matrices.
Given a limit point $p$ for such iterations,
there is a chart in the atlas containing $p$ in its interior,
reducing the study of asymptotic behavior to a matter of local theory.
The construction of the atlas was motivated
by our study of the asymptotics of the Wilkinson shift iteration:
we use bidiagonal coordinates to prove that this well known
algorithm deflates cubically for generic spectra (\cite{LST1})
but only quadratically for certain initial conditions (\cite{LST2}).

Jacobi matrices are frequently parameterized by
its (simple) eigenvalues and the vector $w$
of (positive) first coordinates of its normalized eigenvectors,
the {\it norming constants}.
An algorithm to recover a Jacobi matrix from
these data was known to Stieltjes (\cite{Moser}).
From the procedure, one learns that the set $\JLao \subset \ILa$
of Jacobi matrices with prescribed simple spectrum
is diffeomorphic to $\RR^{n-1}$.
Norming constants break down at the boundary of $\JLao$,
and new techniques are required to study its closure $\JLa$
within the space of symmetric matrices.
In \cite{Tomei}, $\JLa$ was proved to be homeomorphic to
a convex polytope $\PLa$ and the boundary of $\JLa$ was described
as a union of cells of reduced tridiagonal matrices.
But one may go beyond: by making all possible changes
of sign along the off-diagonal entries $(i+1,i)$ of the matrices in $\JLa$,
one obtains $2^{n-1}$ copies of $\JLa$,
which glue along their boundaries to form the compact manifold $\ILa$.

There are significant theoretical advantages for considering
the manifold $\ILa$ instead of $\JLao$ or even $\JLa$.
In algorithms to compute the spectrum of Jacobi matrices,
the limit point is usually a reduced matrix:
if the limit point lies in the interior of the domain,
asymptotic behavior becomes amenable to local theory.
Furthermore, signs of off-diagonal entries are often dropped
in such algorithms.
This procedure, which is computationally practical,
may introduce theoretical complications
akin to inserting absolute values on a smooth function.
Enlarging the domain of such iterations to $\ILa$
may allow for a choice of signs which respects smoothness:
one is then entitled to use Taylor expansions
in the local study of the iteration.

The bidiagonal coordinates 
$\beta^\pi_i$, $i = 1, \ldots, n-1$,
which play the role of generalized norming constants,
are defined on a cover $\UpLa$ of open dense subsets of $\ILa$
indexed by permutations $\pi \in S_n$.
On each $\UpLa$, the bidiagonal coordinates
give rise to a chart of the atlas, 
i.e., a diffeomorphism to $\psi_\pi: \UpLa \to \RR^{n-1}$.
The underlying construction is easy to describe.
A matrix $M$ is \textit{$LU$-positive} if it admits a (unique)
factorization $M = LU$
where $L$ is lower unipotent (i.e., lower triangular with unit diagonal)
and $U$ is upper triangular with positive diagonal entries.
Set $\Lambda = \diag(\lambda_1, \ldots, \lambda_n)$,
$\lambda_1 < \cdots < \lambda_n$ and,
for $\pi \in S_n$, let
$\Lambda^\pi = \diag(\lambda_{\pi(1)}, \ldots, \lambda_{\pi(n)})$.
A matrix $T \in \ILa$ belongs to $\UpLa$ if it admits a diagonalization
$T = Q_\pi^\transpose \Lambda^\pi Q_\pi$
for some orthogonal $LU$-positive matrix $Q_\pi = L_\pi U_\pi$;
in particular, $\Lambda^\pi \in \UpLa$.
Now set $B_\pi = L_\pi^{-1} \Lambda^\pi L_\pi = U_\pi T U_\pi^{-1}$.
From the formulae,
$B_\pi$ is simultaneously lower triangular and upper Hessenberg,
hence lower bidiagonal.
The construction of the chart is complete:
\[ B_\pi = \begin{pmatrix}
\lambda_{\pi(1)} & & & & \\
\beta^\pi_1 & \lambda_{\pi(2)} & & & \\
& \beta^\pi_2 & \lambda_{\pi(3)} & & \\
& & \ddots & \ddots & \\
& & & \beta^\pi_{n-1} & \lambda_{\pi(n)} \end{pmatrix},
\quad \psi_\pi(T) = (\beta^\pi_1, \ldots, \beta^\pi_{n-1}). \]
It turns out that if $T \in \ILa$ is unreduced
then $T \in \UpLa$ for all $\pi \in S_n$ (Lemma \ref{lemma:upla}).

As far as we know, this construction of the matrices $L_\pi$ and $B_\pi$
was first used by Terwilliger in his study
of Leonard pairs (\cite{Terwilliger});
our matrix $L_\pi$, for example, appears in his lemma 4.4 as $E_r$.
Carnicer and Pe\~na (\cite{CP}) also consider changes of basis
leading to bidiagonal matrices in their study of oscillatory matrices.

For Jacobi matrices, bidiagonal coordinates are,
up to a multiplicative factor, quotients of norming constants
(Proposition \ref{prop:nor2bi}):
\[ \beta^\pi_i = \left | \frac{(\lambda_{\pi(i+1)} - \lambda_{\pi(1)})\cdots
(\lambda_{\pi(i+1)} - \lambda_{\pi(i)})}
{(\lambda_{\pi(i)} - \lambda_{\pi(1)})\cdots
(\lambda_{\pi(i)} - \lambda_{\pi(i-1)})} \right |
\frac{w_{\pi(i+1)}}{w_{\pi(i)}},
\quad 1 \le i \le n-1 \]
where $w_{\pi(i)} = |(Q_\pi)_{i,1}|$.
Norming constants, however, yield no chart for a neighborhood
of a diagonal matrix in $\ILa$.
Bidiagonal coordinates imply that in $\UpLa$ appropriate quotients
of norming constants ${w_{\pi(i+1)}}/{w_{\pi(i)}}$
admit natural smooth extensions,
a fact discussed in \cite{Gibson} and \cite{LT}.
There is however no satisfactory definition for the sign of
norming constants for matrices throughout $\ILa$:
this will be discussed more carefully at the end of Section 2.




In the next two sections, we consider the theoretical setup.
Section 2 contains some basic facts about Jacobi matrices
and norming constants,
presented using the concept of $LU$-positivity
so as to prepare the reader to the discussion of bidiagonal coordinates.
We also collect some geometric properties
of the isospectral manifold $\ILa$:
the case $n=3$ is taken as a detailed example.
In section 3, we describe the domains $\UpLa$
both in terms of $LU$-positivity and
based on a cell decomposition of $\ILa$.
We then construct the charts $\psi_\pi: \UpLa \to \RR^{n-1}$
and their inverses $\phi_\pi$:
the bidiagonal coordinates for $T \in \UpLa$
are $(\beta^\pi_1,\ldots,\beta^\pi_{n-1}) = \psi_\pi(T)$.
We also prove that the quotients $\beta^\pi_i/((T)_{i+1,i})$
are smooth strictly positive functions in $\UpLa$.

In order to provide applications, we concentrate on two kinds
of dynamics acting on Jacobi matrices:
\textit{$QR$ steps} (Section 4) and \textit{Toda flows} (Section 5).
Algorithms to compute eigenvalues of Jacobi matrices
which are related to the $QR$ factorization,
as well as the flows in the \textit{Toda hierarchy},
admit a very simple description in bidiagonal coordinates:
they evolve linearly in time.
From this description,
limits at infinity (with asymptotic rates) are immediate.
As a slightly more complicated example,
we prove the cubic convergence of the Rayleigh quotient shift iteration
using a Taylor expansion.
More precisely, given $\Lambda$,
let $G(T) \in \ILa$ be obtained from $T \in \ILa$ by a Rayleigh quotient step:
we prove that there exist $c, \epsilon > 0$
such that if $|(T)_{n,n-1}| < \epsilon$ then
$|(G(T))_{n,n-1}| < \min(\epsilon, c |(T)_{n,n-1}|^3)$.
The reader should compare this argument with
the more complicated study of
the asymptotics of the Wilkinson's shift iteration
in \cite{LST1} and \cite{LST2}.
We conclude the paper with the computation of
the \textit{wave and scattering maps} of the standard Toda flow,
a physical system consisting of $n$ particles
on the line under the influence of a special Hamiltonian.
Moser (\cite{Moser}) had previously computed
the scattering map and Percy Deift (personal communication) the wave map,
but our arguments are significantly different.

We thank the comments presented by the referees,
which led to a much improved text.
The authors gratefully acknowledge support
from CNPq, CAPES, IM-AGIMB and FAPERJ.

\bigbreak

\section{Tridiagonal matrices}

We begin this section by sketching some classical facts about
tridiagonal matrices (\cite{Moser}, \cite{Parlett})
in a phrasing appropriate to our purposes.
Let $\Tt$ be the vector space of real, tridiagonal,
$n \times n$ symmetric matrices.
A matrix $T \in \Tt$ is \textit{Jacobi} (resp. \textit{unreduced})
if $T_{i+1,i} > 0$ (resp. $T_{i+1,i} \ne 0$) for $i = 1, 2, \ldots, n-1$.
Let $\Jj \subset \Tt$ be the open cone of Jacobi matrices and
$\RR^n_{o}$ be the open cone $\{(x_1,\ldots,x_n), x_1 < \cdots < x_n\}$.
The \textit{ordered spectrum map} $\sigma_o$ is defined on the open
set of real symmetric matrices with simple spectrum:
$\sigma_o(S) = (\lambda_1, \ldots, \lambda_n) \in \RR^n_{o}$
lists the eigenvalues of $S$ in increasing order.
Let $O(n)$ be the group of real orthogonal matrices of order $n$.
For an invertible $M$, write the unique {\it $QR$ factorization}
$M = \bQ(M) \bR(M)$,
for $\bQ(M) \in O(n)$ and $\bR(M)$ upper triangular with positive diagonal.
Similarly, when the leading principal minors of $M$ are invertible,
write the {\it $LU$ factorization} $M = \bL(M) \bU(M)$
where $\bL(M)$ is lower unipotent (i.e., lower triangular with unit diagonal)
and $\bU(M)$ is upper triangular.

For a permutation $\pi \in S_n$, consider the matrix $P_\pi$
with $(i,j)$ entry equal $1$ if and only if $i = \pi(j)$
(thus $P_{\pi_1 \pi_2} = P_{\pi_1} P_{\pi_2}$ and $P_\pi e_i = e_{\pi(i)}$).
For $\Lambda = \diag(\lambda_1,\ldots,\lambda_n)$,
with $\lambda_1 < \cdots < \lambda_n$,
set \[ \Lambda^\pi = P_\pi^{-1} \Lambda P_\pi =
\diag(\lambda_{\pi(1)},\lambda_{\pi(2)},\ldots,\lambda_{\pi(n)}) =
\diag(\lambda^\pi_1,\lambda^\pi_2,\ldots,\lambda^\pi_n). \]
Finally, let $\Ee \subset O(n)$ be the group of \textit{sign diagonals}, i.e.,
matrices of the form $E = \diag(\pm 1, \pm 1, \ldots, \pm 1)$.

\begin{defin}
\label{defin:lupositive}
A square matrix $M$ is {\em $LU$-positive}
if $\bU(M)$ is well defined and the diagonal entries of $\bU(M)$ are positive.
Given $T \in \Tt$, $\Lambda = \diag(\sigma_o(T))$ and a permutation $\pi$,
the factorization $T = Q_\pi^\transpose \Lambda^\pi Q_\pi$
is a {\em $\pi$-normalized diagonalization}
if the orthogonal matrix $Q_\pi$ is $LU$-positive.
\end{defin}

Equivalently, $M$ is $LU$-positive if
the determinants of its leading principal minors are positive.
The $\pi$-normalized diagonalization of $T \in \Tt$ is unique
if it exists. Indeed, two factorizations
$Q_1^\transpose \Lambda^\pi Q_1 = Q_2^\transpose \Lambda^\pi Q_2$
yield $Q_2 = E Q_1$ for some $E \in \Ee$,
$Q_1$ and $Q_2$ both $LU$-positive: we must have $E = I$.

We recast a standard result for our purposes.


\begin{theo}
\label{theo:inversevariables}
The eigenvalues of a Jacobi matrix $J$ are distinct.
Given $J$ and a permutation $\pi \in S_n$,
$J$ admits a (unique) $\pi$-normalized diagonalization
$J = Q_\pi^\transpose \Lambda^\pi Q_\pi$
where $\Lambda = \diag(\sigma_o(J))$.
The coordinates of $Q_\pi e_1$ are nonzero.
The map below is a diffeomorphism:
\begin{align*}
\Gamma_\pi: \Jj &\to 
\RR^n_{o} \times
\{ w \in \RR^n \;|\; ||w|| = 1, w_i > 0 \} \\
J &\mapsto \left(\sigma_o(J),
(|(Q_\pi)_{11}|, |(Q_\pi)_{21}|, \ldots, |(Q_\pi)_{n1}|)\right).
\end{align*}
\end{theo}

The second entry $w$ of $\Gamma_\pi(J)$
lists the {\it norming constants} of $J$.
For different $\pi \in S_n$, the coordinates of $w$ are merely permuted.
Indeed, if $J = Q^\transpose \Lambda Q$ and
$J = Q_\pi^\transpose \Lambda^\pi Q_\pi$ then $Q_\pi = EP_\pi^{-1}Q$
for some $E \in \Ee$.

\smallskip

{\nobf Proof:}
Simplicity of the spectrum of $J$ and the fact that the first coordinate
of each eigenvector is nonzero are standard facts (\cite{Parlett}).
In other words, given a diagonalization
$J = {\tilde Q}^\transpose \Lambda^\pi \tilde Q$,
the coordinates of $\tilde w = \tilde Q e_1$ are nonzero.
Clearly, the matrices $\hat Q \in O(n)$ for which
$J = {\hat Q}^\transpose \Lambda^\pi \hat Q$
are of the form $\hat Q = E \tilde Q$ for a sign diagonal $E \in \Ee$,
i.e., we may change signs of rows of $\tilde Q$.
The values of $|(\hat Q)_{i1}|$ do not depend on the choice of $\hat Q$
thus allowing us to define the smooth map $\Gamma_\pi$
using any diagonalization (not necessarily $\pi$-normalized).
We show that one such matrix $\hat Q$ is $LU$-positive
and that $\Gamma_\pi$ is a diffeomorphism by constructing
the inverse of $\Gamma_\pi$.

Construct a Vandermonde matrix $V$ with $V_{ij} = \lambda_{\pi(i)}^{j-1}$ and
a positive diagonal matrix $\tilde W = \diag(\tilde w_1,\ldots,\tilde w_n)$.
The well known formula for the determinant of a Vandermonde matrix
implies that the leading principal minors of $V$ are nonzero.
Thus, there exists a unique sign diagonal $E \in \Ee$
such that $E V$ is $LU$-positive:
the matrices $E \tilde W V = \tilde W E V$ and $\tilde Q = \bQ(E\tilde WV)$
are therefore also $LU$-positive.
We claim that $\tilde J = {\tilde Q}^\transpose \Lambda^\pi \tilde Q$
is a Jacobi matrix.

To prove tridiagonality, we show that
$\tilde J_{ij} = \langle \tilde Je_j, e_i \rangle =
\langle \Lambda^\pi\tilde Q e_j, \tilde Q e_i \rangle$ equals $0$ for $i > j+1$.
For $j = 1, \ldots, n$,
consider the columns $u_j = (\Lambda^\pi)^{j-1}E\tilde w$ and $\tilde q_j$
of $E \tilde WV$ and $\tilde Q$, respectively.
The Krylov subspace $K_j$ spanned by $u_1, \ldots, u_j$ 
is also spanned by the orthonormal vectors $\tilde q_1, \ldots, \tilde q_j$
since $\tilde Q \tilde R = E \tilde WV$
where $\tilde R = \bR(E\tilde WV)$ is upper triangular with positive diagonal.
We have $\Lambda^\pi K_j \subset K_{j+1}$ and therefore 
$\Lambda^\pi \tilde q_j$ is a linear combination of
$\tilde q_1, \ldots, \tilde q_{j+1}$
and $\langle \Lambda^\pi \tilde q_j, \tilde q_i \rangle = 0$ as needed.

We now show that $\tilde J_{j+1,j} > 0$.
The factorization $\tilde Q \tilde R = E\tilde WV$
implies that $\tilde q_i - (\tilde R_{ii})^{-1} u_i \in K_{i-1}$
and similarly $u_{i+1} - \tilde R_{i+1,i+1} \tilde q_{i+1} \in K_i$.
Applying $\Lambda^\pi$ to the first relation we have
$\Lambda^\pi\tilde q_i - (\tilde R_{ii})^{-1} u_{i+1} \in K_i$
and using the second relation we obtain
$\Lambda^\pi\tilde q_i - (\tilde R_{ii})^{-1} \tilde R_{i+1,i+1} \tilde q_{i+1}
\in K_i$.
Now, since $K_i \perp \tilde q_{i+1}$, we have
$\tilde J_{j+1,j} = \langle \Lambda^\pi\tilde q_i, \tilde q_{i+1} \rangle = 
(\tilde R_{ii})^{-1} \tilde R_{i+1,i+1}
\langle \tilde q_{i+1}, \tilde q_{i+1} \rangle
=  (\tilde R_{ii})^{-1} \tilde R_{i+1,i+1} > 0$.

Adding up, $\tilde J = {\tilde Q}^\transpose \Lambda^\pi \tilde Q$ 
is the $\pi$-normalized diagonalization of the Jacobi matrix $\tilde J$.
From this construction,  the inverse of $\Gamma_\pi$ is smooth,
completing the proof.
\qed

\smallskip

Let $\JLao \subset \Jj$ be the set of Jacobi matrices $J$
with $\diag(\sigma_o(J)) = \Lambda$:
norming constants (or, more precisely,
the second coordinate $w$ of $\Gamma_\pi$)
obtain a diffeomorphism between $\JLao$ and
the positive orthant of the unit sphere
$\{ w \in \RR^n \;|\; ||w||=1, w_i >0\}$.
Let $\JLa$ be the closure of $\JLao$ in $\Tt$:
clearly, the boundary of $\JLa$ consists of reduced tridiagonal matrices
with non-negative off-diagonal entries,
including the $n!$ diagonal matrices obtained by permuting
the eigenvalues $\lambda_1, \ldots, \lambda_n$.
The boundary of $\JLa$ is not a smooth manifold:
it has a polytope-like cell structure which
was described in \cite{Tomei}. 
We now give a more explicit description.



For $T \in \JLa$, write a diagonalization
$T = Q^\transpose\Lambda Q$, $Q \in O(n)$.
Consider the matrix $\hat T = Q\Lambda Q^\transpose$:
this matrix is not well defined (due to the sign ambiguity in $Q$)
but the diagonal of $\hat T$ is.
Define $\iota(T)$ to be the diagonal matrix
coinciding with $\hat T$ on the diagonal.
Also, let $\PLa$ be the convex hull of the set of $n!$ matrices $\Lambda^\pi$,
$\pi \in S_n$.


\begin{theo}[\cite{BFR}]
\label{theo:polytope}
The map $\iota$ constructed above
is a homeomorphism $\iota: \JLa \to \PLa$ 
which is a smooth diffeomorphism between interiors.
Furthermore, $\iota(\Lambda^\pi) = \Lambda^{\pi^{-1}}$.
\end{theo}

The original proof of this theorem uses a result of Atiyah
on the convexity of the image
of moment maps defined on K\"ahler manifolds \cite{Atiyah};
a more elementary proof is given in \cite{LT}.
We find this sequence of results to be a good example of the interplay
between high and low roads in linear algebra, so eloquently described
in \cite{Thompson}.

\smallskip

\begin{figure}[ht]
\centerline{\epsfig{height=38mm,file=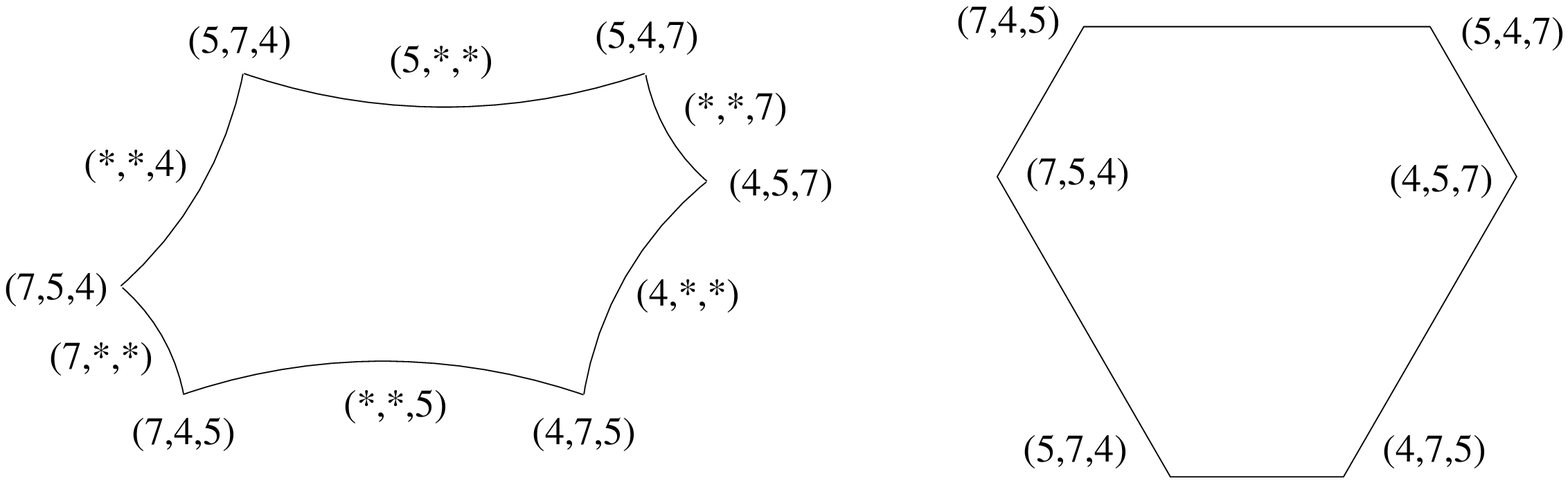}}
\caption{\capsize The spaces $\JLa$ and $\PLa$ for $\Lambda = \diag(4,5,7)$}
\label{fig:hexa}
\end{figure}

We now present in detail the case $n=3$,
where $\JLa$ and $\PLa$ have dimension $2$.
Take $\Lambda = \diag(4,5,7)$;
$\PLa$ is a hexagon contained in the plane $x+y+z = 4+5+7$.
The spaces $\JLa$ and $\PLa$ are given in Figure \ref{fig:hexa}.
The polygon $\PLa$ is drawn in scale;
the drawing of $\JLa$ is schematic.

The triples in the diagram of $\JLa$ are of two kinds.
The vertices, which are diagonal matrices,
are labelled by the three diagonal entries.
The edges have stars in the place of a $2\times 2$ block.
Thus, for instance, the edge $(\ast,\ast,7)$ consists of matrices
of the form
\[ \frac{1}{2} \begin{pmatrix} 9 - \cos 2\theta & \sin 2\theta & 0 \\
\sin 2\theta & 9 + \cos 2\theta & 0 \\ 0 & 0 & 14 \end{pmatrix} =
\begin{pmatrix}c & s & 0 \\ -s & c & 0 \\ 0 & 0 & 1 \end{pmatrix}
\begin{pmatrix} 4 & 0 & 0 \\ 0 & 5 & 0 \\ 0 & 0 & 7 \end{pmatrix}
\begin{pmatrix}c & -s & 0 \\ s & c & 0 \\ 0 & 0 & 1 \end{pmatrix}
\]
where $\theta$ goes from $0$ to $\pi/2$,
$c = \cos\theta$ and $s = \sin\theta$.
Notice that the vertices in $\JLa$ and $\PLa$ have different adjacencies,
in accordance with Theorem \ref{theo:polytope}:
$\Lambda^{\pi_1}$ and $\Lambda^{\pi_2}$ are adjacent in $\PLa$
if and only if $\Lambda^{\pi_1^{-1}}$ and $\Lambda^{\pi_2^{-1}}$
are adjacent in $\JLa$.


\bigskip

Allowing arbitrary signs at off-diagonal entries,
we consider the {\it tridiagonal isospectral manifold} $\ILa$,
the set of real symmetric tridiagonal matrices $T$ which
are conjugate to $\Lambda$.
Define the \textit{sign sequence} of an unreduced matrix $T$
as $\signseq(T) = (\sign(T_{21}), \sign(T_{32}), \ldots, \sign(T_{n,n-1}))$.
The subset of $\ILa$ of unreduced matrices splits into $2^{n-1}$
connected components according to the sign sequence.
Conjugation by sign diagonals takes one component to another.
Thus, $\ILa$ is obtained by gluing $2^{n-1}$ copies of $\JLa$
along their boundaries.

It is shown in \cite{Tomei} 
that $\ILa$ is a compact orientable manifold by proving that
simple spectra are regular values of the restriction to $\Tt$
of the ordered spectrum map $\sigma_o$.

\begin{figure}[ht]
\centerline{\epsfig{height=65mm,file=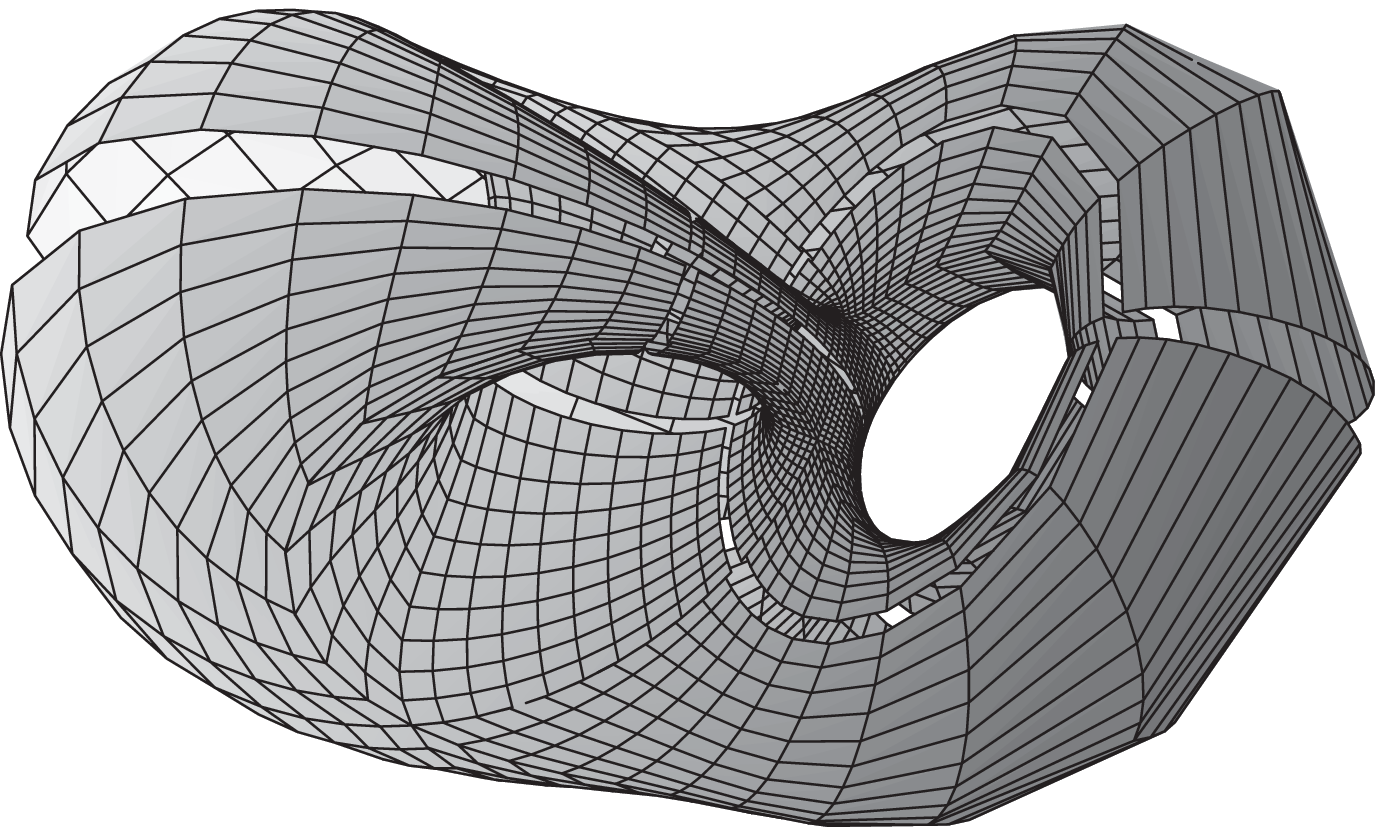}}
\caption{\capsize A 3d rendition of $\ILa$ for $\Lambda = \diag(4,5,7)$}
\label{fig:bitorus}
\end{figure}

For $\Lambda = \diag(4,5,7)$,
Figure \ref{fig:bitorus} shows the manifold $\ILa$, a bitorus.
The vector space $\Tt$ receives an Euclidean metric via the inner product
$\langle T_1, T_2 \rangle = \trace(T_1 T_2)$.
The manifold $\ILa$ is then contained in the intersection
of the hyperplane of matrices of trace $4+5+7$
and the sphere of matrices $T$ with $\langle T, T\rangle = 4^2 + 5^2 + 7^2$:
this intersection is isometric to a sphere centered at the origin in $\RR^4$.
A stereographic projection takes this sphere
(and its subset $\ILa$) to $\RR^3$:
Figure \ref{fig:bitorus} is a snapshot of the image of $\ILa$ under
this projection.
The small gaps were artificially introduced:
these tubular neighborhoods of circles in $\RR^3$
split the manifold into the four hexagons $E\JLa E$, $E \in \Ee$.

\begin{figure}[ht]
\centerline{\epsfig{height=50mm,file=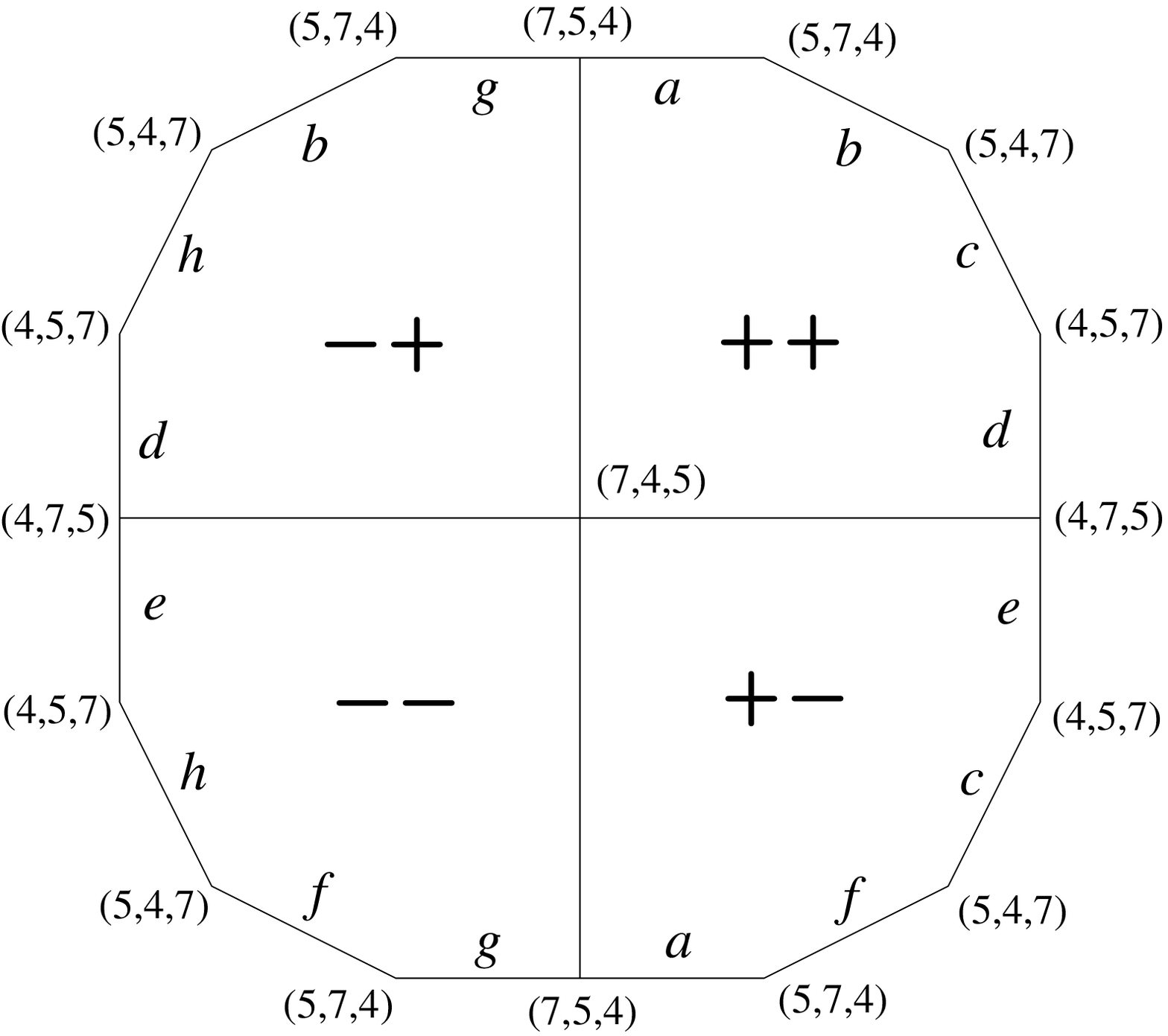}}
\caption{\capsize Gluing instructions for the manifold $\ILa$ }
\label{fig:ilailah}
\end{figure}

Figure \ref{fig:ilailah} shows again $\ILa$ for this example
in a more schematic fashion.
The four hexagons stand for the components
of the subset of unreduced matrices:
sign sequences label the hexagons.
Vertices are diagonal matrices and edges with the same label
are identified.


The bitorus is decomposed as a disjoint union of four open hexagons,
six open edges and six vertices.
We generalize this cell decomposition.
Any tridiagonal matrix $T \in \ILa$ splits into unreduced blocks
$T_1, \ldots, T_k$ along the diagonal.
Consider the subspectra $\Lambda_i = \diag(\sigma_o(T_i))$
and the sign sequences $\signseq(T_i)$.
The \textit{(open) cell containing} $T$ is the subset of $\ILa$
of matrices with the same block partition, subspectra
and sign sequences as $T$.
The cell containing $T$ is naturally identified with
$\Jj_{\Lambda_1} \times \cdots \times \Jj_{\Lambda_k}$
and therefore diffeomorphic to $\RR^{n-k}$,
from Theorem \ref{theo:inversevariables}.
The vertices (or cells of dimension $0$) of $\ILa$
are the diagonal matrices $\Lambda^\pi$
and the set $\JLao$ of Jacobi matrices is the cell of maximal dimension $n-1$
defined by $\signseq(T) = (+,+,\ldots,+)$.

The reader should notice that norming constants do not admit
a smooth natural extension to $\ILa$.
Indeed, again in Figure \ref{fig:ilailah},
norming constants are positive on the $++$ cell of Jacobi matrices.
Crossing the horizontal axis to the cell $+-$ takes the norming
constant $w_2$ (associated with $\lambda_2 = 5$)
through $0$ so, to guarantee smoothness,
we should have signs $+-+$ for the norming constants in the cell $+-$.
On the other hand, we can also go from $++$ to $+-$
through the edges $a$ or $c$ and such crossings
would induce the sign patterns $-++$ and $++-$, respectively.

\bigbreak

\section{Bidiagonal coordinates}

In this section, we construct an atlas for $\ILa$ given by
a family of charts $\psi_\pi$ indexed by the permutations $\pi \in S_n$.
Each chart is a diffeomorphism $\psi_\pi: \UpLa \to \RR^{n-1}$
and the chart domains $\UpLa$,
centered at the diagonal matrices $\Lambda^\pi$
(in the sense that $\psi_\pi(\Lambda^\pi) = 0$),
form an open cover of $\ILa$.
The {\it bidiagonal coordinates} $\beta^\pi_i$ for a matrix $T \in \UpLa$
are the entries of the vector $\psi_\pi(T)$.

For $n = 3$, $\Lambda = \diag(4,5,7)$ and $\pi$ given by
$\pi(1) = 3$, $\pi(2) = 1$, $\pi(3) = 2$ so that $\Lambda^\pi = \diag(7,4,5)$
the set $\UpLa$ is the interior of the polygon in Figure \ref{fig:ilailah}.
Bidiagonal coordinates were used to produce Figure \ref{fig:bitorus}:
$\ILa$ was partitioned into six quadrilaterals centered at diagonal matrices.
This decomposition of $\ILa$ has four vertices (one in the interior of
each of the four hexagons described in the previous section);
each quadrilateral touches each vertex once.
In the figure, the small gaps split each quadrilateral into four smaller ones;
boundaries between quadrilaterals are visible as the lines along which
the mesh loses smoothness.
Lines in each quadrilateral are
level curves of bidiagonal coordinates $\beta^\pi_i$.
We first define the chart domains $\UpLa$.

\begin{defin}
\label{defin:upla}
For a permutation $\pi$,
the {\em chart domain $\UpLa$} is the set of matrices $T \in \ILa$
admitting a $\pi$-normalized diagonalization, i.e.,
the matrices $T$ for which
there exists an $LU$-positive matrix $Q_\pi \in O(n)$
with $T = Q_\pi^\transpose \Lambda^\pi Q_\pi$.
\end{defin}

We now present some properties of the sets $\UpLa$.

\begin{lemma}
\label{lemma:upla}
\begin{enumerate}[(a)]
\item{The sets $\UpLa \subset \ILa$ form an open cover of $\ILa$.}
\item{If $E \in \Ee$ and $T \in \UpLa$ then $ETE \in \UpLa$.}
\item{If $T \in \ILa$ is unreduced
then $T \in \UpLa$ for all permutations $\pi$.
In particular, each set $\UpLa$ is dense in $\ILa$.}
\item{Take $T \in \ILa$ with unreduced blocks $T_1,\ldots,T_k$ of dimensions
$n_1,\ldots,n_k$ along the diagonal.
For a permutation $\pi$, split $\Lambda^\pi$ in blocks:
\[ (\Lambda^\pi)_i = \diag( \lambda^\pi_{n_1+\cdots+n_{i-1}+1}, \ldots,
\lambda^\pi_{n_1+\cdots+n_{i-1}+n_i} ). \]
Then $T \in \UpLa$ if and only if $T_i$ is conjugate to $(\Lambda^\pi)_i$
for $i = 1, \ldots, k$.}
\end{enumerate}
\end{lemma}

\proof
(a) For $T \in \ILa$, write $T = Q^\transpose \Lambda Q$ for some $Q \in O(n)$.
Write the $PLU$ factorization of $Q$, i.e.,
$Q = PLU$ where $P$ is a permutation matrix,
$L$ is lower unipotent and $U$ is upper triangular.
Notice that this is usually possible for $P = P_\pi$
for several permutations $\pi \in S_n$.
Thus, all the leading principal minors of $P_\pi^{-1} Q$ are invertible
and there exists $E \in \Ee$ such that
$Q_\pi = E P_\pi^{-1} Q$ is $LU$-positive
and $T = Q^\transpose \Lambda Q = Q_\pi^\transpose \Lambda^\pi Q_\pi$
belongs to $\UpLa$.
The set of $LU$-positive matrices is open in $\RR^{n \times n}$
and therefore each $\UpLa$ is also open in $\ILa$.

\noindent (b) If $T = Q_\pi^\transpose \Lambda^\pi Q_\pi$ is
the $\pi$-normalized factorization of $T \in \UpLa$ then
$ETE = (EQ_\pi E)^\transpose \Lambda^\pi (EQ_\pi E)$.
The matrix $EQ_\pi E$ is $LU$-positive and therefore $ETE \in \UpLa$.

\noindent (c) The case where $T$ is a Jacobi matrix is discussed
in Theorem \ref{theo:inversevariables}.
If $T$ is unreduced then there exists $E \in \Ee$
such that $ETE$ is Jacobi and item (b) completes the argument.


\noindent (d) Consider a permutation $\pi$,
$T \in \ILa$ with unreduced blocks $T_1,\ldots,T_k$
and the diagonal blocks $(\Lambda^\pi)_i$ as above.
From item (c),
if $T_i$ and $(\Lambda^\pi)_i$ are conjugate then
there exist $LU$-positive matrices $Q_i \in O(n_i)$
with $T_i = Q_i^\transpose (\Lambda^\pi)_i Q_i$.
Let $\tilde Q$ be the matrix with blocks $Q_i$:
the matrix $\tilde Q$ is $LU$-positive and orthogonal;
$T = (\tilde Q)^\transpose \Lambda^\pi \tilde Q$
implies $T \in \UpLa$.

Conversely, assume that $T \in \UpLa$ admits a block decomposition. Then
$T = (P_\pi LU)^{-1} \Lambda (P_\pi LU)$ yields
$UTU^{-1} = L^{-1} \Lambda^\pi L$ and we therefore have
$U_i T_i U_i^{-1} = L_i^{-1} (\Lambda^\pi)_i L_i$
where $U_i$ and $L_i$ are blocks along the diagonal for $U$ and $L$.
\qed

Item (d) yields an alternative description of the chart domain $\UpLa$:
it is the union of all cells in $\ILa$ whose closure contains $\Lambda^\pi$.



We now construct charts $\psi_\pi: \UpLa \to \RR^{n-1}$
and their inverses $\phi_\pi$.

\begin{defin}
\label{defin:phipi}
Given $\beta^\pi_1,\ldots, \beta^\pi_{n-1}$ and
$\Lambda^\pi = \diag(\lambda^\pi_1,\ldots, \lambda^\pi_n)$
build the bidiagonal matrix 
\[ B_\pi = \begin{pmatrix}
\lambda^\pi_1 & & & & \\
\beta^\pi_1 & \lambda^\pi_2 & & & \\
& \beta^\pi_2 & \lambda^\pi_3 & & \\
& & \ddots & \ddots & \\
& & & \beta^\pi_{n-1} & \lambda^\pi_n \end{pmatrix}. \]
Take the diagonalization $B_\pi = L_\pi^{-1} \Lambda^\pi L_\pi$
where $L_\pi$ is lower triangular with unit diagonal.
Define the {\em inverse chart} $\phi_\pi: \RR^{n-1} \to \UpLa$ by
\[ \phi_\pi(\beta^\pi_1,\ldots,\beta^\pi_{n-1}) =
\bQ(L_\pi)^\transpose\,\Lambda^\pi\,\bQ(L_\pi). \]
\end{defin}

Notice that $\phi_\pi(0) = \Lambda^\pi$.
It is easy to see that a matrix $L_\pi$ as above exists;
an explicit formula for its entries is given in the proof
of Proposition \ref{prop:nor2bi}.
The claim that $\UpLa$ is a valid counterdomain
for $\phi_\pi$ requires a proof.
Set $Q_\pi = \bQ(L_\pi)$:
the matrix $T = \phi_\pi(\beta^\pi_1,\ldots,\beta^\pi_{n-1}) =
Q_\pi^\transpose\,\Lambda^\pi\,Q_\pi$ is clearly symmetric.
On the other hand, the factorization
$T = \bR(L_\pi) B_\pi (\bR(L_\pi))^{-1}$ implies that $T$ is upper Hessenberg
and therefore $T \in \ILa$.
Since $\bR(L_\pi)$ has positive diagonal and
$Q_\pi = L_\pi (\bR(L_\pi))^{-1}$,
we have that the matrix $Q_\pi$ is $LU$-positive,
$T = Q_\pi^\transpose\,\Lambda^\pi\,Q_\pi$ is
a $\pi$-normalized diagonalization and therefore $T \in \UpLa$.

\begin{defin}
\label{defin:psipi}
For $T \in \UpLa$, take its $\pi$-normalized diagonalization
$T = Q_\pi^\transpose \Lambda^\pi Q_\pi$ and
write $L_\pi = \bL(Q_\pi)$, $U_\pi = \bU(Q_\pi)$, $R_\pi = U_\pi^{-1}$. Set
\[ B_\pi = R_\pi^{-1} T R_\pi = L_\pi^{-1} \Lambda^\pi L_\pi. \]
The matrix $B_\pi$ is bidiagonal
and its off-diagonal entries $\beta^\pi_i = (B_\pi)_{i+1,i}$,
$i = 1, \ldots, n-1$,
are the {\em $\pi$-bidiagonal coordinates} of $T$.
The {\em chart} $\psi_\pi: \UpLa \to \RR^{n-1}$
is the map taking $T$ to $(\beta^\pi_1,\ldots,\beta^\pi_{n-1})$.
\end{defin}

We prove that $B_\pi$ is indeed bidiagonal.
From $B_\pi = R_\pi^{-1} T R_\pi$, $B_\pi$ is upper Hessenberg and
from $B_\pi = L_\pi^{-1} \Lambda^\pi L_\pi$, it is lower triangular
with diagonal entries $\lambda^\pi_1, \ldots, \lambda^\pi_n$.


The maps $\psi_\pi$ and $\phi_\pi$ are clearly smooth.
By construction, one is the inverse of the other,
implying the following result.

\begin{theo}
\label{theo:bidiagonal}
The map $\phi_\pi: \RR^{n-1} \to \UpLa$
is a diffeomorphism with inverse $\psi_\pi: \UpLa \to \RR^{n-1}$.
\end{theo}

\smallskip

As an example of bidiagonal coordinates, let $\Lambda = \diag(4,5,7)$.
Set $\pi(1) = 3$, $\pi(2) = 1$, $\pi(3) = 2$.
Matrices will be described by their $\pi$-bidiagonal coordinates
$x = \beta_1^\pi$ and $y = \beta_2^\pi$.
Since $B_\pi = L_\pi^{-1} \Lambda^\pi L_\pi$, we obtain
\[
\Lambda^\pi = \begin{pmatrix} 7 & 0 & 0 \\ 0 & 4 & 0 \\
0 & 0 & 5 \end{pmatrix}, \quad
B_\pi = \begin{pmatrix} 7 & 0 & 0 \\ x & 4 & 0 \\
0 & y & 5 \end{pmatrix}, \quad
L_\pi = \begin{pmatrix} 1 & 0 & 0 \\ -x/3 & 1 & 0 \\
- xy/2 & y & 1 \end{pmatrix} \]
and writing $Q_\pi = L_\pi U_\pi$ we have
\[ Q_\pi =
\frac{1}{r_1r_2}
\begin{pmatrix}
{6}{r_2} & {6x(2+3y^2)} & {xy}{r_1} \\
{-2x}{r_2} & {3(12+x^2y^2)} & {-6y}{r_1} \\
{-3xy}{r_2} & {2y(18-x^2)} & {6}{r_1}
\end{pmatrix}, \]
where $r_1 = \sqrt{36 + 4x^2 + 9x^2y^2}$
and $r_2 = \sqrt{36 + 36y^2 + x^2y^2}$.
From $T = Q_\pi^\transpose \Lambda^\pi Q_\pi$, we have
$T = \Lambda^\pi + M/(r_1^2r_2^2)$ where
\[ M =
\begin{pmatrix}
-6x^2(2 + 3y^2) r_2^2 & 6x r_2^3 & 0 \\
6x r_2^3 &
(72 + 108y^2) r_1^2 - (72 - 4x^2) r_2^2 &
6y r_1^3 \\
0 & 6y r_1^3 & -2y^2(18-x^2) r_1^2
\end{pmatrix}. \]
This is an explicit parametrization $\phi_\pi: \RR^2 \to \UpLa \subset \ILa$
of the polygon in Figure \ref{fig:ilailah}.
From this formula, $x = 0$ implies $(T)_{11} = 7$ and $(T)_{21} = 0$
while $y = 0$ gives $(T)_{33} = 5$, $(T)_{32} = 0$,
consistent with the description of $\UpLa$ in Lemma \ref{lemma:upla}.

\smallskip

\begin{prop}
\label{prop:nor2bi}
For any permutation $\pi$ and any Jacobi matrix $J \in \JLao$,
the norming constants $w_i$ and the
$\pi$-bidiagonal coordinates $\beta^\pi_i$ are related by
\[ w_{\pi(i)} = w_{\pi(1)} \left | \frac{\beta^\pi_1 \cdots \beta^\pi_{i-1}}%
{(\lambda^\pi_i - \lambda^\pi_1) \cdots (\lambda^\pi_i - \lambda^\pi_{i-1})}
\right |,
\quad 2 \le i \le n, \]
\[ \beta^\pi_i = \left | \frac{(\lambda_{\pi(i+1)} - \lambda_{\pi(1)})\cdots
(\lambda_{\pi(i+1)} - \lambda_{\pi(i)})}
{(\lambda_{\pi(i)} - \lambda_{\pi(1)})\cdots
(\lambda_{\pi(i)} - \lambda_{\pi(i-1)})} \right |
\frac{w_{\pi(i+1)}}{w_{\pi(i)}},
\quad 1 \le i \le n-1. \]
\end{prop}

\proof
Let $\Lambda^\pi$ and $B_\pi$ be as above.
Set
\[ L = \begin{pmatrix}
1 & 0 & 0 & \cdots & 0 \\ \\
\frac{\beta^\pi_1}{\lambda^\pi_2 - \lambda^\pi_1} & 1 & 0 & \cdots & 0 \\ \\
\frac{\beta^\pi_1 \beta^\pi_2}%
{(\lambda^\pi_3 - \lambda^\pi_1)(\lambda^\pi_3 - \lambda^\pi_2)} &
\frac{\beta^\pi_2}{\lambda^\pi_3 - \lambda^\pi_2}  & 1 & \cdots & 0 \\
\vdots & \vdots & \vdots & & \vdots \\
\frac{\beta^\pi_1 \beta^\pi_2 \cdots \beta^\pi_{n-1}}%
{(\lambda^\pi_n - \lambda^\pi_1)(\lambda^\pi_n - \lambda^\pi_2)\cdots%
(\lambda^\pi_n - \lambda^\pi_{n-1})} &
\frac{\beta^\pi_2 \cdots \beta^\pi_{n-1}}%
{(\lambda^\pi_n - \lambda^\pi_2)\cdots (\lambda^\pi_n - \lambda^\pi_{n-1})} &
& \cdots & 1
\end{pmatrix}.
\]
A straightforward computation verifies that
$L B_\pi = \Lambda^\pi L$ and therefore $L = L_\pi$.

Norming constants are given by the absolute values
of entries in the first column of $Q_\pi$ which in turn is
the normalization of the first column of $L_\pi$.
\qed



\smallskip

Bidiagonal coordinates change signs together with off-diagonal entries
in a simple fashion.

\begin{lemma}
\label{lemma:ETE}
If $E = \diag(\sigma_1, \ldots, \sigma_n) \in \Ee$, $T \in \UpLa$, and
$ \psi_\pi(T) = (\beta^\pi_1,\ldots,\beta^\pi_{n-1}) $
then
\( \psi_\pi(ETE) =
(\sigma_1 \sigma_2 \beta^\pi_1,\ldots,
\sigma_{n-1} \sigma_n \beta^\pi_{n-1})\).
In other words, if $B_\pi$ and $\tilde B_\pi$ are the bidiagonal matrices
associated to $T$ and $ETE$ (as in Definition \ref{defin:psipi}), respectively,
then $\tilde B_\pi = E B_\pi E$.
\end{lemma}

\proof
From Lemma \ref{lemma:upla}, if $T \in \UpLa$ then $\tilde T = ETE \in \UpLa$.
Clearly, if the $\pi$-normalized decomposition of $T$ is
$T = Q_\pi^\transpose \Lambda^\pi Q_\pi$ then
$ETE = (EQ_\pi E)^\transpose \Lambda^\pi (EQ_\pi E)$ is the
$\pi$-normalized diagonalization of $ETE$.
Also, $\bL(EQ_\pi E) = E \bL(Q_\pi) E$
and therefore, if $L_\pi = \bL(Q_\pi)$ then
$\tilde B_\pi = (E L_\pi E)^{-1} \Lambda^\pi (E L_\pi E) = E B_\pi E$.
\qed

Norming constants break down at the boundary of $\JLa$.
We will prove in the following proposition, however,
that near a reduced matrix $T$ with $(T)_{i+1,i} = 0$
the values of $(T)_{i+1,i}$ and $\beta^\pi_i$ are comparable.
This will be useful when we use bidiagonal coordinates
to study the asymptotics of isospectral maps.


\begin{prop}
\label{prop:bbeta}
Given $\Lambda$ and $\pi$, the quotient $q^\pi_i: \UpLa \to \RR$
defined by $q^\pi_i(T) = \beta^\pi_i(T)/((T)_{i+1,i})$
is smooth, positive and $q^\pi_i(\Lambda^\pi) = 1$.
Also, $q^\pi_i(ETE) = q^\pi_i(T)$ for all $T \in \UpLa$ and $E \in \Ee$.
\end{prop}

In particular $\beta^\pi_i$ and $(T)_{i+1,i}$
have the same sign regardless of $\pi$.

\smallskip

\proof
Clearly, $(T)_{i+1,i} = 0$ if and only if $ETE = T$
where $E_{j,j} = 1$ (resp. $-1$) for $j \le i$ (resp. $j \ge i+1$).
Thus, from Lemma \ref{lemma:ETE},
$\beta^\pi_i(T) = 0$ if and only if $(T)_{i+1,i} = 0$.

We study the $i$-th partial derivative of 
$(\phi_\pi(\beta^\pi_1, \ldots, \beta^\pi_{n-1}))_{i+1,i}$
when $\beta^\pi_i = 0$.
Take a path $T: \RR \to \ILa$,
$T(t) = \phi_\pi(\beta^\pi_1, \ldots, t, \ldots, \beta^\pi_{n-1})$.
From Lemma \ref{lemma:ETE}, all entries of $T(t)$ 
except $(i+1,i)$ (and $(i,i+1)$) are even functions of $t$
and therefore the corresponding entries of $T'(0)$ equal $0$.
On the other hand, since $\phi_\pi$ is a diffeomorphism, 
$T'(0)$ must be nonzero.
It follows that $(T'(0))_{i+1,i} \ne 0$
and therefore $q^\pi_i(T)$ is well defined, smooth
and nonzero even when the denominator vanishes, i.e.,
at reduced matrices.




The symmetry property indicated in the last claim
follows from Lemma \ref{lemma:ETE}.
Positivity is obvious for Jacobi matrices, extends to unreduced
matrices by symmetry and to reduced matrices by continuity.
In order to compute $q^\pi_i(\Lambda^\pi)$,
consider the path $T(t) = \phi_\pi(0,\ldots, t, \ldots, 0)$
(with $t$ in the $i$-th position).
Clearly, $B_\pi(t) = \Lambda^\pi + t (e_{i+1})^\ast e_i$
so that $\beta^\pi_i(t) = t$.
A straightforward computation yields
\[ (T(t))_{i+1,i} = \frac{(\lambda^\pi_{i+1} - \lambda^\pi_i)^2}%
{(\lambda^\pi_{i+1} - \lambda^\pi_i)^2 + t^2}\;t;
\quad
q^\pi_i(T(t)) = 
\frac{(\lambda^\pi_{i+1} - \lambda^\pi_i)^2 + t^2}%
{(\lambda^\pi_{i+1} - \lambda^\pi_i)^2}. \]
The result now follows by setting $t = 0$.
\qed

\smallskip

\section{Iterations in $\UpLa$}

We now apply bidiagonal coordinates to the study of the dynamics
of $QR$ type iterations.
As a simple example, we present in Theorem \ref{theo:rayleigh}
a new proof of the well known fact that the Rayleigh quotient shift iteration
has cubic convergence.
A subtler example is the Wilkinson's shift iteration
which is studied with the same technique in \cite{LST1} and \cite{LST2}.
Excellent references for the spectral theory of Jacobi matrices
are \cite{Demmel}, \cite{Golub} and \cite{Parlett}.

\smallskip


For an open neighborhood $X \subset \RR$
of the spectrum $\{\lambda_1, \ldots, \lambda_n\} = \sigma(\Lambda)$
and a continuous function $f: X \to \RR$ taking nonzero values
on $\sigma(\Lambda)$,
there is a smooth map $F: \ILa \to \ILa$,
the \textit{$QR$ step induced by $f$}, given by
\[ F(T) = \bQ(f(T))^\transpose\,T\,\bQ(f(T)). \]
Continuity of $f$ is sufficient to imply that $F$
is a well defined smooth function:
indeed, if $f$ and the polynomial $p$ coincide on $\sigma(\Lambda)$
then $f(T) = p(T)$ for all $T \in \ILa$ and therefore $F = P$,
the $QR$ step induced by $p$.
The standard $QR$ step corresponds to $f(x) = x$.
Since $T$ and $f(T)$ commute, we also have
\[ F(T) = \bR(f(T))\,T\,\bR(f(T))^{-1}. \]
From the first formula, $F(T)$ is symmetric;
from the second, it is upper Hessenberg with
sub-diagonal elements with the same signs as in $T$.
Thus, $F: \ILa \to \ILa$ preserves $\JLao$, the other cells
and the open subsets $\UpLa$.

Let $F^\be =
\phi^{-1}_\pi \circ F \circ \phi_\pi: \RR^{n-1} \to \RR^{n-1}$;
in other words, $F^\be$ is obtained from $F|_{\UpLa}$ 
by a change of variables using bidiagonal coordinates.


\begin{prop}
\label{prop:stepbi}
For $f$ taking nonzero values on the spectrum of $T$,
\[ F^\be(\beta^\pi_1, \ldots, \beta^\pi_{n-1}) =
\left( \left | \frac{f(\lambda_{\pi(2)})}{f(\lambda_{\pi(1)})} \right |
\beta^\pi_1, \ldots,
\left | \frac{f(\lambda_{\pi(n)})}{f(\lambda_{\pi(n-1)})} \right |
\beta^\pi_{n-1} \right). \]
Also, $QR$ iterations generically converge to diagonal matrices.
More precisely, if $T \in \UpLa$ and
\( |f(\lambda_{\pi(1)})| > |f(\lambda_{\pi(2)})| > \cdots >
|f(\lambda_{\pi(n)})| \),
then $\lim_{k \to +\infty} F^k(T) = \Lambda^\pi$ with asymptotics
\[ \beta^\pi_i = \lim_{k \to +\infty} (F^k(T))_{i+1,i}
\left|\frac{f(\lambda_{\pi(i)})}{f(\lambda_{\pi(i+1)})}\right|^k. \]
\end{prop}

\proof
Take $T \in \UpLa$ and set $T'= F(T) = \bQ(f(T))^\transpose T \bQ(f(T))$.
We show that $T'\in \UpLa$ and relate
the corresponding matrices $B_\pi$ and $B'_\pi$.
Consider the $\pi$-normalized diagonalization
$T = Q_\pi^\transpose \Lambda^\pi Q_\pi$ and write $L_\pi = \bL(Q_\pi)$.
Since $f(T) = Q_\pi^\transpose f(\Lambda^\pi) Q_\pi$
and $\bQ(ZM) = Z \bQ(M)$
for an arbitrary matrix $Z \in O(n)$ and invertible $M$, we have
$T' = (\tilde Q)^\transpose \Lambda^\pi \tilde Q$ where
\( \tilde Q = Q_\pi \bQ(Q_\pi^\transpose f(\Lambda^\pi) Q_\pi)
= \bQ( f(\Lambda^\pi) Q_\pi)\).
Take $Q'_\pi = \bQ( |f|(\Lambda^\pi) Q_\pi)$:
clearly, $Q'_\pi$ is $LU$-positive and 
$Q'_\pi = E\tilde Q$ for some $E \in \Ee$ and therefore
$T' = (Q'_\pi)^\transpose \Lambda^\pi Q'_\pi$
is the $\pi$-normalized diagonalization of $T' \in \UpLa$.
Write $L_\pi = \bL(Q_\pi)$, $L'_\pi = \bL(Q'_\pi)$.
Since $\bL(M) = \bL(M R)$ for arbitrary $LU$-positive matrices $M$
and invertible, upper triangular $R$,
we have $\bL(M) = \bL(\bQ(M))$ and thus
$L'_\pi = \bL( |f|(\Lambda^\pi) Q_\pi)$.
Notice that if $D$ is an invertible diagonal matrix and $M$ is $LU$-positive
then $\bL(DM) = D \bL(M) D^{-1}$:
we obtain $L'_\pi = |f|(\Lambda^\pi) L_\pi (|f|(\Lambda^\pi))^{-1}$
and therefore
%
%
\[ B'_\pi = (L'_\pi)^{-1} \Lambda^\pi L_\pi =
(|f|(\Lambda^\pi)) B_\pi (|f|(\Lambda^\pi))^{-1}. \]
This finishes the proof of the first formula.
The convergence properties now follow easily from Proposition \ref{prop:bbeta}.
\qed

This proposition yields yet another evidence for the naturality
of the bidiagonal coordinates $\beta^\pi_i$.

The cubic convergence to deflation
of the $QR$ iteration with Rayleigh quotient shift 
is well known (\cite{Parlett});
using bidiagonal coordinates, we deduce it from a Taylor expansion.
For $s \in \RR$, let $f_s(x) = x - s$ so that
the $QR$ step $F_s: \ILa \to \ILa$ is defined for $s \notin \sigma(\Lambda)$.
In other words, we have a map
$\bF: (\RR \smallsetminus \sigma(\Lambda)) \times \ILa \to \ILa$,
$\bF(s,T) = F_s(T)$.
The map $\bF$ cannot be continuously extended to $\RR \times \ILa$;
it follows from Proposition \ref{prop:stepbi}, however, that $\bF$
can be continuously extended to pairs $(s,T)$ if $s = \lambda_i$
and the (possibly reduced) matrix $T$ has the eigenvalue $\lambda_i$
in the spectrum of its bottom block.
More formally, consider the set
\[ \Dd_{\bF} = ((\RR \smallsetminus \sigma(\Lambda)) \times \ILa) \cup
\bigcup_{i = 1}^n 
\left( \{\lambda_i\} \times \bigcup_{\pi \in S_n, \pi(i) = n} \UpLa \right)
\subset \RR \times \ILa. \]
The set $\Dd_{\bF}$ is open since
points $T \in \{\lambda_{\pi(n)}\} \times \UpLa$ 
admit the explicit open neighborhood 
$(\lambda^\pi_n - \gamma/2, \lambda^\pi_n + \gamma/2) \times \UpLa \subset \Dd_{\bF}$
where $\gamma = \min_{i \ne j} |\lambda_i - \lambda_j|$
is the spectral gap of $\Lambda$.
The function $\bF$ is defined in $\Dd_{\bF}$ by
\[
\bF(s,\phi_\pi(\beta^\pi_1, \ldots, \beta^\pi_{n-1})) =
\phi_\pi\left( \left|\frac{\lambda^\pi_2 - s}{\lambda^\pi_1 - s}\right| \beta^\pi_1,
\ldots,
\left|\frac{\lambda^\pi_n - s}{\lambda^\pi_{n-1} - s}\right| \beta^\pi_{n-1}
\right). \eqno{(\ast)} \]
The \textit{$QR$ iteration with Rayleigh quotient shift} 
$G: \Dd_G \to \ILa$ is defined by $G(T) = \bF((T)_{n,n},T)$
where $\Dd_G \subseteq \ILa$ is the open set
$\{ T \in \ILa \;|\; ((T)_{n,n},T) \in \Dd_{\bF} \}$.
Notice that, from a numerical point of view,
falling outside $\Dd_G$ is an instant win: $(T)_{n,n}$ is an eigenvalue.

The \emph{deflation set} $\Delta_0 \subset \ILa$ is the set of matrices
$T$ with $(T)_{n,n-1} = 0$.
The set $\Delta_0$ is the disjoint union of the subsets $\Delta_0^i$
of matrices $T \in \Delta_0$ with $(T)_{n,n} = \lambda_i$.
Notice that $\Delta_0^i$ is diffeomorphic to $\Tt_{\Lambda_{\hat i}}$
where $\Lambda_{\hat i} =
\diag(\lambda_1, \ldots, \lambda_{i-1}, \lambda_{i+1}, \ldots, \lambda_n)$
and therefore a connected component of $\Delta_0$.
Clearly, $\Delta_0 \subset \Dd_G$
since $T_{n,n-1} = 0$ and $T_{n,n} = \lambda^\pi_n$
imply that $T \in \UpLa$.
The \emph{(compact) deflation neighborhood} $\Delta_\epsilon$
is the set of matrices $T \in \ILa$ with $|(T)_{n,n-1}| \le \epsilon$.

\begin{theo}
\label{theo:rayleigh}
There exist $\epsilon > 0$ and $c > 0$
such that $\Delta_\epsilon \subset \Dd_G$,
$G(\Delta_\epsilon) \subset \Delta_\epsilon$ and
for all $T \in \Delta_\epsilon$ we have
$ |(G(T))_{n,n-1}| \le c\;|(T)_{n,n-1}|^3 $.
\end{theo}

\proof
Since the chart domains $\UpLa$ are open dense sets covering $\ILa$,
there exist compact sets $K_\pi \subset \UpLa$ with
$\bigcup_{\pi \in S_n} K_\pi = \ILa$.
Take $M_\pi > 0$ such that $\psi_\pi(K_\pi)$
is contained in the box $(-M_\pi,M_\pi)^{n-1}$:
we then have 
\[ \Delta_0 = \bigcup_{\pi \in S_n} \phi_\pi((-M_\pi, M_\pi)^{n-2} \times \{0\}). \]
Equation $(\ast)$ above yields a formula for
$G^{\phi_\pi}: \Dd_{G^{\phi_\pi}} \to \RR^{n-1}$
where $\Dd_{G^{\phi_\pi}} = \psi_\pi(\Dd_G \cap \UpLa) \subseteq \RR^{n-1}$
is an open set with $\RR^{n-2} \times \{0\} \subset \Dd_{G^{\phi_\pi}}$ and 
\[ s = s(\beta^\pi_1, \ldots, \beta^\pi_{n-1}) =
(\phi_\pi(\beta^\pi_1, \ldots, \beta^\pi_{n-1}))_{n,n}. \]
For each $\pi$ there exist $\epsilon'_\pi > 0$ such that,
if $A'_\pi = (-M_\pi, M_\pi)^{n-2} \times (-\epsilon'_\pi,\epsilon'_\pi)$
then $(\beta^\pi_1, \ldots, \beta^\pi_{n-1}) \in A'_\pi$
implies $|s - \lambda^\pi_n| < \gamma/2$
(where $\gamma$ is the spectral gap)
and therefore $\overline{A'_\pi} \subset \Dd_{G^{\phi_\pi}}$.
Due to the presence of absolute values, the function $G^{\phi_\pi}$
is almost certainly not smooth in $A'$.
Set $g_\pi: A'_\pi \to \RR$,
\[ g_\pi(\beta^\pi_1, \ldots, \beta^\pi_{n-1}) = 
\frac{\lambda^\pi_n - s}{|\lambda^\pi_{n-1} - s|}. \]
The function $g_\pi$ is smooth
and $\beta^\pi_{n-1} = 0$ implies $g_\pi(\beta^\pi_1, \ldots, \beta^\pi_{n-1}) = 0$.
Since $g$ is even (from Lemma \ref{lemma:ETE}),
its first order partial derivatives at
points with $\beta^\pi_{n-1} = 0$ all vanish
and the Taylor expansion for $g_\pi$ at such points
starts with terms of degree $2$.
By compactness of $\overline{A'_\pi}$,
there exists a constant $c_\pi$ such that
\[g_\pi(\beta^\pi_1, \ldots, \beta^\pi_{n-1}) \le c_\pi (\beta^\pi_{n-1})^2, \quad
|(G^{\phi_\pi}(\beta^\pi_1, \ldots, \beta^\pi_{n-1}))_{n-1}| \le
c_\pi |\beta^\pi_{n-1}|^3. \]
From Proposition \ref{prop:bbeta}, $\beta^\pi_{n-1}$ and $(T)_{n,n-1}$
are comparable: there exists $\tilde c_\pi$
such that $T \in \phi_\pi(A'_\pi)$ implies
\[ |(G(T))_{n,n-1}| \le \tilde c_\pi |T_{n,n-1}|^3. \]
Take $c = \max_{\pi \in S_n} \tilde c_\pi$ and
$\epsilon > 0$ such that $\epsilon < c^{-1/2}$ and
$\Delta_\epsilon \subset A'$.
If $T \in \Delta_\epsilon$ we therefore have
$|(G(T))_{n,n-1}| \le c |(T)_{n,n-1}|^3 \le \epsilon$,
proving the claims.
\qed

\bigskip

\section{Toda flows}

Recall that the Toda flow
(\cite{Flaschka}, \cite{Moser}, \cite{DNT})
solves the differential equation
\[ J'(t) = [ J(t) , \Pi_a (J(t))], \quad J(0) = J_0. \]
Here the bracket is the usual Lie bracket
on matrices $[A_1,A_2] = A_1A_2- A_2A_1$ and
$\Pi_a(M)$ is the skew-symmetric matrix having the same lower
triangular entries as $M$. 
As is well known, this flow preserves spectrum and the set $\JLa$.
If $w(t)$ is the vector of norming constants for $J(t)$, we have
\[ w(t) = \frac{\exp(t\Lambda) w(0)}{||\exp(t\Lambda) w(0)||}: \]
thus, up to normalization, the function $w$ is the solution
of a linear differential equation.
Taking quotients and using Proposition \ref{prop:nor2bi}
shows that the evolution of $\beta^\pi_i$ is truly linear:
\[ \frac{d}{dt}\beta^\pi_i(t) =
(\lambda_{\pi(i+1)} - \lambda_{\pi(i)}) \beta^\pi_i(t). \]
In other words, $B_\pi' = [B_\pi,-\Lambda^\pi]$.

Clearly, the Toda flow is well defined in $\ILa$.
Similar formulae hold for other flows in the Toda hierarchy:
for a function $g$, consider the differential equation
$T' = [T, \Pi_a g(T)]$;
it turns out that,
despite $g(T)$ not being tridiagonal,
$[T, \Pi_a g(T)]$ is symmetric and tridiagonal.
Integrate the differential equation to define $T(t)$ for all $t \in \RR$.
In bidiagonal coordinates, it is easy to compute limits and asymptotics
of Toda flows.

\begin{prop}
\label{prop:todabi}
In $\pi$-bidiagonal coordinates,
the equation $T' = [T, \Pi_a g(T)]$
becomes a decoupled linear system: \( \frac{d}{dt}\beta^\pi_i =
(g(\lambda^\pi_{i+1}) - g(\lambda^\pi_i)) \beta^\pi_i\).
In particular, if $T(0) \in \UpLa$ and
\( g(\lambda_{\pi(1)}) > g(\lambda_{\pi(2)}) > \cdots
> g(\lambda_{\pi(n)}) \)
then 
\( \lim_{t \to +\infty} \beta^\pi_k(t) = 0 \),
so that $\lim_{t \to +\infty} T(t) = \Lambda^\pi$
with asymptotics
\[ \beta^\pi_k(0) = \lim_{t \to +\infty} (T(t))_{k,k+1}
\exp\left(\left(g(\lambda_{\pi(k)}) - g(\lambda_{\pi(k+1)})\right) t \right).\]
\end{prop}


\proof
The formula below follows by direct computation
(\cite{Kostant}, \cite{LT}, \cite{Symes}, \cite{Symes2}):
\[ T(t) = \bQ(\exp(t\,g(T(0))))^\transpose\,T(0)\,\bQ(\exp(t\,g(T(0)))), \]
or, in other words, $\bT(t) = F(T(0))$ where $f(x) = \exp(t\,g(x))$.
From Proposition \ref{prop:stepbi},
\[ (\beta^\pi_1(t), \ldots, \beta^\pi_{n-1}(t)) =
\left(
e^{(g(\lambda^\pi_2) - g(\lambda^\pi_1)) t} \beta^\pi_1(0), \ldots,
e^{(g(\lambda^\pi_n) - g(\lambda^\pi_{n-1})) t} \beta^\pi_{n-1}(0) \right) \]
and the differential equation for $\beta^\pi_i$ follows by taking derivatives.
The last formula is now a consequence of Proposition \ref{prop:bbeta}.
\qed






As an application of the bidiagonal variables we consider the
scattering properties of the Toda flow.
From a more physical point of view,
the Toda flow is the evolution of $n$
particles of mass $1$ on the line given by the Hamiltonian
\[ H(x,y) =
\sum_{k=1}^n \frac{y_k^2}{2}  + \sum_{k=1}^{n-1} \exp( x_k - x_{k+1} ), \]
where $x=(x_1,\ldots,x_n)$ and $y=(y_1,\ldots,y_n)$ are respectively the
positions and velocities of the particles.
Without loss of generality,
\[ \sum_k x_k(t) = \sum_k y_k(t) = 0. \]
More explicitly, positions and
velocities satisfy the differential equation
\[ x_k' = H_{y_k} = y_k, \quad
y_k' = - H_{x_k} = \exp( x_{k-1} - x_k) - \exp(x_k - x_{k+1}), \quad
k =1,\ldots,n, \]
where we take the formal boundary conditions
$x_0 = -\infty$, $x_{n+1} = \infty$.
The two versions of the Toda flow are related by Flaschka's transformation
(\cite{Flaschka}):
\[ J_{k,k} = - \frac{1}{2} y_k, \quad
J_{k,k+1} = \frac{1}{2} \exp\left(\frac{x_k - x_{k+1}}{2}\right). \]
Notice that $\trace(J) = 0$ for the Jacobi matrix $J$
constructed from $x$ and $y$.

%

We know that $J(t)$ tends to $\Lambda^{\pi_\pm}$
when $t \to \pm\infty$ where
$\pi^-$ is the identity permutation and 
$\pi^+$ is the reversal $\pi^+(k) = n+1-k$.
The convergence of the diagonal entries implies that
the velocities $y_k(t)$ approach $-2\lambda_{\pi_\pm(k)}$.
From the convergence to $0$ of the off-diagonal entries,
the force $\pm \exp(x_k - x_{k+1})$ between particles $k$ and $k+1$
tends to $0$ when $t \to \pm\infty$.
Thus, asymptotically, the particles undertake independent uniform motions
of the form $x_k(t) \approx c_k^\pm t + d_k^\pm$ when $t \to \pm\infty$.
Clearly, $c_k^\pm = -2\lambda_{\pi_\pm(k)}$.
The \textit{wave} and \textit{scattering maps} are
\begin{align*}
W^\pm(x_1(0),\ldots,x_n(0),y_1(0),\ldots,y_n(0)) &=
(c_1^\pm,\ldots,c_n^\pm,d_1^\pm,\ldots,d_n^\pm), \\
S(c_1^-,\ldots,c_n^-,d_1^-,\ldots,d_n^-) &=
(c_1^+,\ldots,c_n^+,d_1^+,\ldots,d_n^+),
\end{align*}
respectively. They are related by $S \circ W^- = W^+$.
Moser (\cite{Moser}) proved that $S$ is indeed well defined and computed it.
We obtain Moser's result by first computing wave maps
in bidiagonal coordinates.

\begin{prop}
\label{prop:moser}
Given initial conditions $(x_1(0),\ldots,x_n(0),y_1(0),\ldots,y_n(0))$
with $\sum x_k(0) = \sum y_k(0) = 0$, apply Flaschka's transformation
to obtain $J(0)$ with eigenvalues $\lambda_1 < \cdots < \lambda_n$
and bidiagonal coordinates $\beta_k^{\pi_\pm}$ for
$\pi_-(k) = k$ and $\pi_+(k) = n+1-k$.
Then
\[ c_k^\pm = -2\lambda_{\pi_\pm(k)}, \quad
d_k^\pm = \sum_{j < k} \frac{-2j}{n}\tilde\beta_j^{\pi_\pm}
+ \sum_{j \ge k} \frac{2(n-j)}{n}\tilde\beta_j^{\pi_\pm} + (n-2k+1)\log 2, \]
where $\tilde\beta_k^{\pi_\pm} = \log \beta_k^{\pi_\pm}$.
The scattering map is given by 
\[ c_{n+1-k}^+ = c_k^-, \quad
d_{n+1-k}^+ = d_k^-
+ 2 \sum_{j < k} \log|c_j^- - c_k^-|
- 2 \sum_{j > k} \log|c_j^- - c_k^-|. \]
\end{prop}

\proof
From Proposition \ref{prop:todabi} with $g(z) = z$,
\[ J_{k,k+1} = \frac{1}{2} \exp\left(\frac{x_k - x_{k+1}}{2} \right) \approx
\exp\left( (\lambda_{\pi_\pm(k+1)} - \lambda_{\pi_\pm(k)} )t \right)
\beta_k^{\pi_\pm} \]
in the sense that quotients tend to $1$ when $t \to \pm\infty$.
Taking logs,
\[ (x_k - x_{k+1}) - 2 (\lambda_{\pi_\pm(k+1)} - \lambda_{\pi_\pm(k)} )t 
\approx (d_k^\pm - d_{k+1}^\pm) = 2 \tilde\beta_k^{\pi_\pm} + 2 \log 2, \]
where now the difference goes to zero.
Since $\sum_k x_k = 0$, we have $\sum_k d_k^\pm = 0$ and
the formula for the wave operator follows.

Proposition \ref{prop:nor2bi} yields
$\tilde{\beta}_k^+ = - \tilde{\beta}_{n-k}^- + \delta_k$
where
\[ \delta_k =  \left(\sum_{j > n-k} \epsilon_{n-k,j}
- \sum_{j < n-k} \epsilon_{n-k,j}\right)
- \left(\sum_{j > n+1-k} \epsilon_{n+1-k,j}
- \sum_{j < n+1-k} \epsilon_{n+1-k,j}\right) \]
and $\epsilon_{ij} = \log\left|\lambda_{\pi^-(i)}- \lambda_{\pi^-(j)}\right|$.
Use this equation to write
\begin{align*}
d_k^+ &= - \left( \sum_{j < k} \frac{-2j}{n}\tilde{\beta}_{n-j}^-
+ \sum_{j \ge k} \frac{2(n-j)}{n}\tilde{\beta}_{n-j}^- \right) \\
&\qquad + \sum_{j < k} \frac{-2j}{n} \delta_j
+ \sum_{j \ge k} \frac{2(n-j)}{n} \delta_j
+ (n-2k+1)\log 2.
\end{align*}
Replacing $k$ by $r(k)= n+1-k$ in the equation for $d_k^-$ yields
\[ d_{n+1-k}^- = - \left( \sum_{j < k} \frac{-2j}{n}\tilde{\beta}_{n-j}^-
+ \sum_{j \ge k} \frac{2(n-j)}{n}\tilde{\beta}_{n-j}^- \right) \\
- (n-2k+1)\log 2 \]
and therefore
\[ d_k^+ = d_{n+1-k}^- + \sum_{j < k} \frac{-2j}{n} \delta_j
+ \sum_{j \ge k} \frac{2(n-j)}{n} \delta_j
+ 2(n-2k+1)\log 2. \]
The simplification yielding the scattering map is now an easy exercise.
\qed

\bigbreak

\bigskip\bigskip\bigbreak

{

\parindent=0pt
\parskip=0pt
\obeylines

Ricardo S. Leite, Departamento de Matem\'atica, UFES
Av. Fernando Ferrari, 514, Vit\'oria, ES 29075-910, Brazil

\smallskip

Nicolau C. Saldanha and Carlos Tomei, Departamento de Matem\'atica, PUC-Rio
R. Marqu\^es de S. Vicente 225, Rio de Janeiro, RJ 22453-900, Brazil

\smallskip

rsleite@pq.cnpq.br
nicolau@mat.puc-rio.br; http://www.mat.puc-rio.br/$\sim$nicolau/
tomei@mat.puc-rio.br

}

\end{document}